\documentclass[11pt]{amsart}%
\usepackage{amsfonts}
\usepackage{amsmath}
\usepackage{amssymb}
\usepackage{graphicx}%
\newtheorem{theorem}{Theorem}[section]
\theoremstyle{plain}

\numberwithin{equation}{section}
\begin{document}
\title[Homogeneous Solutions]{Homogeneous Solutions to Fully Nonlinear Elliptic Equations}
\author{Nikolai NADIRASHVILI}
\address{Department of Mathematics\\
University of Chicago\\
5734 S. University Ave., Chicago, IL 60637}
\curraddr{LATP\\
Centre de Math\'{e}matiques et Informatique \\
39, rue F. Joliot-Curie \\
13453 Marseille Cedex \\
France }
\email{nicholas@math.uchicago.edu}
\author{Yu YUAN}
\address{Department of Mathematics, Box 354350\\
University of Washington\\
Seattle, WA 98195}
\email{yuan@math.washington.edu}
\thanks{Both authors are partially supported by NSF grants, and the second author also
by a Sloan Research Fellowship}

\begin{abstract}
We classify homogeneous degree $d\neq2$ solutions to fully nonlinear elliptic equations.

\end{abstract}
\maketitle

In this note, we show that any homogeneous degree other than $2$ solution to
fully nonlinear elliptic equations must be \textquotedblleft harmonic".
Consider the fully nonlinear elliptic equation $F\left(  D^{2}u\right)  =0$
with $\mu I\leq\left(  F_{ij}\right)  =\left(  F_{M_{ij}}\left(  M\right)
\right)  \leq\mu^{-1}I.$ Nirenberg [N] derived the a priori $C^{2,\alpha}$
estimates for the above equation in dimension $2$ in 1950s. Krylov [K] and
Evans [E] showed the same a priori estimates for the above equations in
general dimensions under the assumption that $F$ is convex. As a modest
investigation of a priori estimates for general fully nonlinear elliptic
equations without convexity condition, we study the homogeneous solutions.

\begin{theorem}
Let $u$ be a continuous in $\mathbb{R}^{n}\backslash\left\{  0\right\}  $
homogeneous degree $d\neq2$ solution to the elliptic equation $F\left(
D^{2}u\right)  =0$ in $\mathbb{R}^{n}$ with $F\in C^{1}.$ Then $u$ is harmonic
in a possible new coordinate system in $R^{n},$ namely%
\[
\sum_{i,j=1}^{n}F_{ij}\left(  0\right)  D_{ij}u\left(  x\right)  =0.
\]
Consequently, $u\equiv0$ if $-\left(  n-2\right)  \,<d<0$ or $d$ is not an
integer; otherwise $u$ is a homogeneous harmonic polynomial with integer
degree $d.$
\end{theorem}

In contrast to the variational problem, Sverak and Yan [SY] constructed
homogeneous degree less than $1$ minimizers to some strongly convex
functional. Also Safonov [S] constructed homogeneous order $\alpha$ $\in(0,1)$
solutions to linear non-divergence elliptic equations with variable
coefficients earlier on.

As one simple application to special Lagrangian equations [HL] $F\left(
D^{2}u\right)  =\sum_{i=1}^{n}\arctan\lambda_{i}-c=0,$ where $\lambda_{i}s$
are the eigenvalues of the Hessian $D^{2}u.$ It follows from our theorem that
any homogeneous degree other than $2$ solutions must be a harmonic polynomial
(and it also forces $c=0$).

When $d\in\lbrack0,1+\alpha\left(  n,\mu\right)  ),$ our theorem follows from
Krylov-Safonov $C^{\alpha}$ estimates (cf. [CC, corollary 5.7]). The missing
case $d=2$ is delicate. One only knows that any homogeneous degree $2$
solution to the above fully nonlinear elliptic equation in dimension $3$ is
quadratic [HNY, p. 426].

Now we show our theorem.

\begin{proof}
We first consider the case that $u$ is smooth in $\mathbb{R}^{n}%
\backslash\left\{  0\right\}  .$ Set $\sum=\left\{  \left\vert x\right\vert
^{d-2}D^{2}u\left(  \frac{x}{\left\vert x\right\vert }\right)  |x\in
\mathbb{R}^{n}\backslash\left\{  0\right\}  \right\}  $ and $\Gamma=\left\{
M|F\left(  M\right)  =0\right\}  .$ For the homogeneous order $d$ function
$u\left(  x\right)  ,$ $D^{2}u\left(  x\right)  =\left\vert x\right\vert
^{d-2}D^{2}u\left(  \frac{x}{\left\vert x\right\vert }\right)  .$ Let
$\left\vert x\right\vert \rightarrow0$ for $d>2$ or $\left\vert x\right\vert
\rightarrow\infty$ for $d<2,$ we see that $0\in\sum.$ Also $u$ is a solution
to $F\left(  D^{2}u\right)  =0,$ then the cone $\sum\subseteq\Gamma.$ Now
$F\in C^{1}$ and $\left(  F_{ij}\left(  0\right)  \right)  >0,$ we know that
the unique tangent plane of $\Gamma$ at $0$ includes $\sum.$ It follows that
$\sum\perp\left(  F_{ij}\left(  0\right)  \right)  ,$ or%
\[
\sum_{i,j=1}^{n}F_{ij}\left(  0\right)  D_{ij}u\left(  x\right)  =0.
\]
Without loss of generality, we assume $\left(  F_{ij}\left(  0\right)
\right)  =I$ through out the proof, then
\[
0=\bigtriangleup u\left(  x\right)  =\left\vert x\right\vert ^{d-2}\left[
d\left(  d+n-2\right)  u\left(  \frac{x}{\left\vert x\right\vert }\right)
+\bigtriangleup_{S^{n-1}}u\left(  \frac{x}{\left\vert x\right\vert }\right)
\right]  .
\]
The remaining conclusion of the theorem follows.

Next we show the regularity of the viscosity solution $u$ away from $0.$ Set
$\lambda=d\left(  d+n-2\right)  $ and $\theta=x/\left\vert x\right\vert .$ To
start, we prove that $u\left(  \theta\right)  $ is a viscosity solution to%
\begin{equation}
\bigtriangleup_{S^{n-1}}u+\lambda u=0. \label{Eeigen}%
\end{equation}
Let any smooth $\varphi\left(  \theta\right)  $ touch $u$ from the above at
$\theta_{0},$%
\begin{align*}
\varphi &  \geq u\ \ \text{in a neighborhood of }\theta_{0}\\
\varphi\left(  \theta_{0}\right)   &  =u\left(  \theta_{0}\right)  ,
\end{align*}
then%
\begin{align*}
\left\vert x\right\vert ^{d}\varphi\left(  \frac{x}{\left\vert x\right\vert
}\right)   &  \geq\left\vert x\right\vert ^{d}u\left(  \frac{x}{\left\vert
x\right\vert }\right)  \ \ \text{in a neighborhood of }\theta_{0}\\
\left\vert x\right\vert ^{d}\varphi\left(  \theta_{0}\right)   &  =\left\vert
x\right\vert ^{d}u\left(  \theta_{0}\right)  .
\end{align*}
From our assumption that $u$ is a viscosity (sub) solution, it follows that%
\[
F\left(  D^{2}\left(  \left\vert x\right\vert ^{d}\varphi\left(  \frac
{x}{\left\vert x\right\vert }\right)  \right)  \right)  \geq0
\]
or%
\[
F\left(  \left\vert x\right\vert ^{d-2}D_{x}^{2}\varphi\left(  \theta\right)
\right)  \geq0.
\]

Let $\left\vert x\right\vert \rightarrow0$ for $d>2$ or $\left\vert
x\right\vert \rightarrow\infty$ for $d<2,$ we see that $F\left(  0\right)
\geq0.$ If we use the fact $u$ is also a viscosity (super) solution, we can
derive that $F\left(  0\right)  \leq0.$ So $F\left(  0\right)  =0,$ and
\[
\frac{F\left(  tD_{x}^{2}\varphi\left(  \theta\right)  \right)  -F\left(
0\right)  }{t}\geq0.
\]
Let $t\rightarrow0,$ we see that%
\[
\sum_{i,j=1}^{n}F_{ij}\left(  0\right)  D_{ij}\varphi\left(  \frac
{x}{\left\vert x\right\vert }\right)  \geq0
\]
or%
\[
\bigtriangleup_{S^{n-1}}\varphi+\lambda\varphi\geq0.
\]
Thus $u$ is a viscosity sub solution to (\ref{Eeigen}). Similarly, $u$ is a
viscosity super solution to the same equation.

Let $N_{\varepsilon}$ be an $\varepsilon$ neighborhood of any $\theta_{0}$ on
$S^{n-1},$ with $\varepsilon$ small enough so that $N_{\varepsilon}$ is in a
narrow strip, then there exists positive smooth function $h$ on
$N_{\varepsilon}$ such that%
\[
\bigtriangleup_{S^{n-1}}h+\lambda h\leq0.
\]
Let $\psi$ be the smooth solution to (\ref{Eeigen}) in $N_{\varepsilon}$ with
the boundary value $u$ on $\partial N_{\varepsilon},$ then $q=\frac{\psi-u}%
{h}$ is a viscosity solution to%
\[
\bigtriangleup_{S^{n-1}}q+2\frac{\nabla h}{h}\cdot\nabla q+\frac
{\bigtriangleup_{S^{n-1}}h+\lambda h}{h}q=0,
\]
where $\nabla h\cdot\nabla q$ simply denotes some linear combinations of first
order derivatives of $q$ in some local coordinates for $N_{\varepsilon}$,
which we avoid for the sake of simple notation. Now that the coefficient
$\frac{\bigtriangleup_{S^{n-1}}h+\lambda h}{h}\leq0,$ it follows from
[W,Corollary 3.20] that $q=0$ in $N_{\varepsilon}.$ Therefore, $u$ is smooth
in $N_{\varepsilon}$ and then on the whole $S^{n-1}.$
\end{proof}

\end{document}